\input amstex 
\documentstyle{amsppt}
\input bull-ppt
\keyedby{bull272/cxv}

\define\bbc{\Bbb C}
\define\bbd{\Bbb D}
\define\bbn{\Bbb N}
\define\sco{\scr O}
\define\op{\operatorname}
\define\scx{\bold X}
\define\scz{\bold Z}
\define\lb{\lambda}

\topmatter
\cvol{26}
\cvolyear{1992}
\cmonth{April}
\cyear{1992}
\cvolno{2}
\cpgs{276-279}
\title Analytic varieties versus integral varieties\\
of Lie algebras of vector fields\endtitle
\author Herwig Hauser and Gerd M\" uller\endauthor
\shorttitle{Analytic varieties versus integral varieties}
\address Mathematisches Institut der Universit\" at 
Innsbruck, A-6020 
Innsbruck, Austria\endaddress
\address Fachbereich Mathematik der Universit\" at Mainz, 
D-6500 Mainz, Germany\endaddress
\date April 16, 1991 and, in revised form, August 28, 
1991\enddate
\thanks This work was done during a visit of the second 
author at the
University of Innsbruck. He thanks the members of the 
Mathematics Department
for their hospitality\endthanks
\abstract We associate to any germ of an analytic variety 
a Lie algebra of
tangent vector fields, the {\it tangent algebra}. 
Conversely, to any Lie
algebra of vector fields an analytic germ can be 
associated, the {\it integral
variety}. The paper investigates properties of this 
correspondence: The set of
all tangent algebras is characterized in purely Lie 
algebra theoretic terms.
And it is shown that the tangent algebra determines the 
analytic type of the
variety.\endabstract
\endtopmatter

\document

Local analytic varieties, defined as zero sets of complex 
analytic functions,
can equally be considered as integral varieties associated 
to certain Lie
algebras of vector fields. This is the theme of the 
present note. As a
consequence one obtains a new way of studying 
singularities of varieties by
looking at their Lie algebra. It turns out that the Lie 
algebra determines
completely the variety up to isomorphism. Thus one may 
replace, to a certain
extent, the local ring of functions on the variety by the 
Lie algebra of
vector fields tangent to the variety.

We shall give a brief account of these observations. 
Details will appear
elsewhere, see \cite{HM1, HM2}. The paper of Omori \cite 
O, which treats the
same topic in a special case, served us as a valuable 
source of inspiration.
Various ideas are already apparent there.

Consider a germ $X$ of a complex analytic variety embedded 
in some smooth
ambient space, $X\subset(\bbc^n,0)$. In this note, germ of 
variety shall
always mean reduced but possibly reducible complex space 
germ. We associate to
$X$ the Lie algebra $\bbd_X$ of vector fields on 
$(\bbc^n,0)$ tangent to $X$.
To do so let $\bbd$ denote the Lie algebra of germs of 
analytic vector fields
on $(\bbc^n,0)$. We identify $\bbd$ with 
$\op{Der}\,\sco_n$, the Lie algebra of
derivations of the algebra $\sco_n$ of germs of analytic 
functions
$(\bbc^n,0)\to\bbc$. We then set
$$
\bbd_X=\{D\in\bbd,\ D(I_X)\subset I_X\},
$$
where $I_X\subset\sco_n$ is the ideal of functions 
vanishing on $X$. This is a 
subalgebra of $\bbd$. It will be called the {\it tangent 
algebra\/} of $X$. In case $X$
is a nonreduced germ, simple examples show that the 
tangent algebra of $X$ and
of its reduction $X_{\op{red}}$ may coincide. This limits 
our interest to the
reduced case. In this context, two main problems arise:

$\bullet$ Characterize all Lie subalgebras $A\subset\bbd$ 
that are of the form
$A=\bbd_X$ for a suitable $X$.

$\bullet$ Find out to what extent the abstract Lie algebra 
$\bbd_X$ determines
the variety $X$. 

\heading 1. Tangent algebras were characterized by Lie 
algebra properties\endheading

In order to discuss the first problem let us fix some 
notation. A subalgebra
$A$ of a Lie algebra $B$ will be called {\it balanced} (in 
$B$) if $A$ contains
no ideal $\ne 0$ of $B$ but an element $a\ne 0$ such that
$$
[a,B]\subset A\quad\text{ and }\quad[[a,B],B]\subset A.
$$
A {\it visible\/} subalgebra of $B$ is a subalgebra $A$ 
that admits a chain
of subalgebras
$$
A=A_m\subset A_{m-1}\subset\cdots\subset A_0=B
$$
such that $A_k$ is maximal balanced in $A_{k-1}$ for 
$k=1,\dots,m$. In case
$m=1$, i.e., if $A$ is a maximal balanced subalgebra of 
$B$, $A$ is called 
{\it maximal visible}. Note that these notions are of a 
purely Lie algebra
theoretic nature.

For a finite family $\scx=\{X_1,\dots,X_p\}$ of germs 
$X_i\subset(\bbc^n,0)$
let $\bbd_{\scx}=\bigcap_i\bbd_{X_i}$ be the Lie algebra 
of vector fields
tangent to all $X_i$ (the $X_i$ may be contained in each 
other). Our first
result may be considered as a variation of the classical 
Frobenius Theorem in
the singular case (see e.g., \cite{N, 2.11}).

\thm{Theorem 1} Let $A\subset\bbd$ be a subalgebra.

{\rm (a)} There is a set of germs $\scx$ as above such 
that $A=\bbd_{\scx}$ if
and only if $A$ is a visible subalgebra of $\bbd$.

{\rm (b)} There is a smooth germ $X\subset(\bbc^n,0)$ 
different from
$\emptyset$ and $(\bbc^n,0)$ such that $A=\bbd_X$ if and 
only if $A$ is a
maximal visible subalgebra of $\bbd$.

{\rm (c)} There is an irreducible germ 
$X\subset(\bbc^n,0)$ with an isolated
singularity at $0$ such that $A=\bbd_X$ if and only if $A$ 
is a maximal
visible subalgebra of the algebra $\bbd_0$ of vector 
fields vanishing at $0$.

{\rm(d)} There is an analytic germ $X\subset(\Bbb C^n,0)$ 
such that
$A=\Bbb D_X$ if and only if $A$ is geometric in $\Bbb D$, 
i.e. by
definition, $A$ is visible in every subalgebra $B$ of 
$\Bbb D$
containing $A$.
\ethm

\rem{Comments} (i) It is easy to see that the family 
$\scx$ of germs $X_i$ is
not unique. For example, if $\scx$ is the set of 
irreducible components $X_i$
of some germ $X=\bigcup X_i$ one has $\bbd_X=\bbd_{\scx}$. 
Moreover,
$\bbd_X=\bbd_{X,\op{Sing}\, X}$ where $\op{Sing}\,X$ 
denotes the singular
subspace of $X$. But in case $\scx$ is an irredundant set 
of irreducible germs,
i.e., deleting any germ from $\scx$ alters $\bbd_{\scx}$, 
the family $\scx$ is
uniquely determined by $\bbd_{\scx}$. In particular, the 
variety $X$
of a maximal geometric subalgebra as in (b) and (c) is 
unique.

(ii) There is a relative version of Theorem 1 where $\bbd$ 
is replaced by
$\bbd_{\scz}$ for some set of germs $\scz$ and where all 
varieties
$\scx$ associated to visible subalgebras of $\bbd_{\scz}$ 
are determined. Namely, a
subalgebra $A$ of $\bbd_{\scz}$ is visible if and only if 
there is a set of
irreducible germs $\scx$ with $X_i\nsubset Z_j$ for all 
$i,j$ such that
$A=\bbd_{\scz,\scx}$. Theorem 1 represents the cases 
$\scz=\emptyset$, resp\.
$\scz=\{0\}$. 
\endrem

\heading 2. Singularities are determined by their tangent 
algebra\endheading

We now turn to the second problem, the characterization of 
the isomorphism type
of germs via their Lie algebra. If $X,Y\subset(\bbc^n,0)$ 
are isomorphic then
the associated Lie algebras $\bbd_X$ and $\bbd_Y$ are 
isomorphic. In fact,
every isomorphism $X\to Y$ can be extended to\ an 
automorphism $\phi$ of
$(\bbc^n,0)$ with algebra automorphism 
$\phi^*\:\sco_n\to\sco_n$. Then
$\Phi(D):=\phi^*\circ D\circ(\phi^*)^{-1}$ defines an 
automorphism $\Phi$ of
$\bbd$ with $\Phi(\bbd_Y)=\bbd_X$. By abuse of notation we 
write again
$\Phi=\phi^*$. This map is continuous if $\bbd$ is 
provided with the topology
induced from the coefficientwise topology on $\sco_n$. 
Conversely we have

\thm{Theorem 2} Let $X$ and $Y$ be germs of analytic 
varieties in $(\bbc^n,0)$
different from $\emptyset$. Assume that $n\ge 3$. For 
every isomorphism
$\Phi\:\bbd_Y\to\bbd_X$ of topological Lie algebras there 
is a unique
automorphism $\phi$ of $(\bbc^n,0)$ sending $X$ onto $Y$ 
and such that
$\Phi=\phi^*$.
\ethm

Thus the analytic isomorphism type of $X$ is entirely 
given by the abstract
topological Lie algebra $\bbd_X$. Omori \cite O proved 
this in
the special case of
weighted homogeneous varieties. 

We indicate some ideas appearing in the proofs of Theorems 
1 and 2.

\heading 3. Proof of Theorem 1\endheading 

In order to study visible subalgebras of $\bbd$ we 
associate to any
$A\subset\bbd$ the germ $X(A)$ in $(\bbc^n,0)$ defined by 
the ideal
$\sqrt{I(A)}$ of $\sco_n$ where
$$
I(A)=\{g\in\sco_n,\ g\cdot\bbd\subset A\}.
$$
Here the $\sco_n$-module structure of $\bbd$ is used. The 
germ $X(A)$ will be
called the {\it integral variety\/} of $X$. Note that 
every germ
$X\subset(\bbc^n,0)$ different from $\emptyset$ and 
$(\bbc^n,0)$ can be
recovered from $\bbd_X$ as $X=X(\bbd_X)$: The inclusion 
$X(\bbd_X)\subset X$ is
obvious from the definition. For the converse, assume that 
some $g\in
I(\bbd_X)$ does not belong to $I_X$. Consider the vector 
fields
$g\partial_{x_1},\dots,g\partial_{x_n}$. In every point 
outside the zero set of $g$
in $X$ they are linearly independent. As they are tangent 
to $X$ by definition
of $I(\bbd_X)$ a Theorem of Rossi \cite{R, Theorem 3.2} 
implies that the germ
of $X$ taken in such a point is isomorphic to 
$(\bbc^n,0)$. But these points
are dense in $X$ and we get a contradiction.

Let us now consider assertion (b) of Theorem 1. The proof 
that $\bbd_X$ is a 
balanced subalgebra of $\bbd$ is a bit involved and will 
be left out.
Concerning maximality, assume that $\bbd_X$ is contained 
in a balanced
subalgebra $A\subset\bbd$. Then in fact $\bbd_X\subset 
A\subset\bbd_{X(A)}$.
One shows that $A$ balanced implies $X(A)\ne\emptyset$ and 
$(\bbc^n,0)$.
Moreover $X(A)=X(\bbd_{X(A)})\subset X(\bbd_X)=X$. Now if 
$X$ is smooth one
deduces from $\bbd_X\subset\bbd_{X(A)}$ that $X(A)=X$. 
This shows $\bbd_X=A$
and proves necessity in (b).

For sufficiency, start with a maximal visible subalgebra 
$A\subset\bbd$.
Similarly as above $A\subset\bbd_{X(A)}$ with 
$X(A)\ne\emptyset$ and
$(\bbc^n,0)$. As $\bbd_{X(A)}$ is balanced, maximality of 
$A$ gives
$A=\bbd_{X(A)}$. Write $X=X(A)$. If 
$\op{Sing}\,X\ne\emptyset$ then
$\bbd_{\op{Sing}\,X}$ is balanced. Again by maximality, 
the inclusion
$\bbd_X\subset\bbd_{\op{Sing}\,X}$ is actually an 
equality. This implies
$X=\op{Sing}\,X$, which is impossible. Therefore $X$ is 
smooth.

Part (a) of Theorem 1 is proved by induction. Here one 
proves and uses at once
the relative version of the theorem mentioned earlier. To 
illustrate, let
$\scx=\{X\}$ consist of one singular germ $X$. Choose 
$k\in\bbn$ maximal with
$Z=\op{Sing}^kX:=\op{Sing}(\cdots(\op{Sing}(X))\ne%
\emptyset$. The inclusion
$\bbd_X\subset\bbd$ is split into 
$\bbd_X=\bbd_{Z,X}\subset\bbd_Z$ and
$\bbd_Z\subset\bbd$. The first is visible by induction and 
the second is
maximal visible by (b) since $Z$ is smooth.

Conversely, if $A$ is a visible subalgebra of $\bbd$ use 
induction on the
length of the chain and the relative version of part (b) 
to find $\scx$.

\heading 4. Proof of Theorem 2\endheading 

We conclude with some remarks on the proof of Theorem 2. 
For $f\in\sco_n$
consider the $\bbc$-linear map $\lb_f\:\bbd_X\to\bbd_X$ 
defined by
$$
\lb_f(D)=\Phi(f\cdot\Phi^{-1}(D)).
$$
If $\Phi\:\bbd_Y\to\bbd_X$ is induced from an automorphism 
$\phi$ of
$(\bbc^n,0)$, say $\Phi=\phi^*$, one checks by computation 
that the equality
$\lb_f(D)=\phi^*(f\,)\cdot D$ holds for all $D$ in 
$\bbd_X$. If $\Phi$ is an
arbitrary continuous Lie algebra isomorphism, we are led 
to establish the same
equality in order to recover a map $\phi$ that could be an 
appropriate
candidate to induce $\Phi$ and to define an isomorphism 
between $X$ and $Y$.

Thus the first thing to do is to check whether any vector 
field $D$ is mapped
by $\lb_f$ into the $\sco_n$-module $(D)$ generated by 
$D$. This can be seen
for all $D$ of a certain dense subset $U$ of 
$I_X\cdot\bbd$ by writing $(D)$ as
an intersection of subalgebras of $\bbd_X$ of form 
$\bbd_{X,Z}$. This is the key
step in the proof and it is here that we need the 
assumption $n\ge 3$. Once
this is accomplished, the relative version of Theorem 1 
and the fact that
$\Phi$ is a Lie algebra isomorphism guarantee that $\lb_f$ 
maps 
$\bbd_{X,Z}$ into
$\bbd_{X,Z}$. Hence the module $(D)$ is mapped into 
itself. This implies
$$
\lb_f(D)=\phi^*(f,D)\cdot D
$$
with suitable factor $\phi^*(f,D)\in\sco_n$. Then the 
continuity of $\Phi$ is
used to show that $\phi^*(f,D)$ is actually independent of 
$D$, say
$\phi^*(f,D)=\phi^*(f\,)$. Therefore, again by continuity,
$$
\lb_f(D)=\phi^*(f\,)\cdot D
$$
will hold for all $D\in I_X\cdot\bbd$. Finally we deduce 
from this equality
that the map $\phi$ thus obtained is an automorphism of 
$(\bbc^n,0)$ mapping
$X$ to $Y$ and inducing $\Phi$.

\enddocument